\documentclass{amsart}

\usepackage{a4,xy,amssymb}
\usepackage{amsthm,amsmath,amscd,amssymb,graphics}
\usepackage[hyperindex,pdfmark]{hyperref}
\xyoption{all}

\makeindex

\def\cR{{\mathcal R}}

\def\C{{\mathbb C}}

\def\Z{{\mathbb Z}}
\def\N{{\mathbb N}}

\def\eop{\hspace*{\fill}$\square$}
\def\implies{\Rightarrow}

\def\to{\rightarrow}

\def\<{\langle}
\def\>{\rangle}

\theoremstyle{definition}
\newtheorem{defn}{Definition}[section]
\newtheorem{thm}{Theorem}[section]
\newtheorem*{thm*}{Theorem}

\newtheorem{lemma}[thm]{Lemma}
\theoremstyle{remark}

\newtheorem*{prf}{Proof}
\newcounter{kd}

\newcommand{\grp}[1]{\ensuremath{\mathsf{#1}}}

\newcommand{\Rep}{\ensuremath{\mathsf{Rep}}}
\newcommand{\Ker}{\ensuremath{\mathsf{Ker}}}
\newcommand{\Dim}{\ensuremath{\mathsf{Dim}}}
\newcommand{\iss}{\ensuremath{\mathsf{iss}}}
\newcommand{\GL}{\ensuremath{\mathsf{GL}}}
\newcommand{\Repr}{\ensuremath{\mathsf{Rep}}}

\newcommand{\Mat}{\ensuremath{\mathsf{Mat}}}

\newcommand{\se}[1]{\begin{equation*}\begin{split}#1\end{split}\end{equation*}}

\newcommand {\vtx}[1]{*+[o][F-]{\scriptstyle{#1}}}
\newcommand {\vierkant}[1]{*+[F-]{\scriptstyle{#1}}}

\newcommand {\unit}{\mathsf {1\!\!\!1}}
\newcommand {\Tr}{\mathsf {Tr}}

\author{Raf Bocklandt}
\address{Universiteit Antwerpen (UIA) \\ B-2610 Antwerp (Belgium)}
\email{rbockl@uia.ac.be}
\begin{document}

\bibliographystyle{plain}
\title{Smooth quiver quotient varieties}
\begin{abstract}
In this paper we classify all the quivers and corresponding dimension vectors having
a smooth space of semisimple representation classes.
The result is that these quiver settings can be reduced via some specific reduction steps
to 3 simple types.
\end{abstract}

\maketitle
\section{Introduction and motivation}

Many problems in representation theory
can be reduced to representations of quivers.
Suppose $A$ is a finitely generated algebra and
$\Repr_n A$ is the space of $n$-dimensional complex representations of $A$.
On this space is an action of $\GL_\alpha$ and one can divide out this action
by taking the affine quotient
to obtain a new space $\iss_n A:=\Repr_n A/\GL_\alpha$
classifying the equivalence classes of $n$-dimensional semisimple representations of $A$.
(see \cite{Kraft:1984})

If $W\in \Rep_n A$ is a semisimple representation, one can wonder
what the structure of $\iss_n A$ around the point $p$
corresponding to the equivalence class of $W$ looks like. If $W$
is a smooth point in $\Rep_n A$ there is a neighborhood of $p$
that is \'etale (or analytically) isomorphic to a neighborhood of
the zero representation in the quotient space, $\iss_{\alpha_p},
Q_p$ of a quiver setting $(Q_p,\alpha_p)$ which is called the
local quiver setting of $p$. This local quiver setting depends on
the structure of $W$ as a direct sum of simple representations
\[
W := S_1^{\oplus a_1}\oplus \cdots \oplus S_k^{\oplus a_k}.
\]
(For the exact construction see \cite{LBP})

So if one for example asks whether $\iss_n A$ is smooth in the point $p$ one can as well
ask whether its local quiver setting has a quotient space that is smooth in zero.
As we will see below this is the same as asking whether this quotient space is an affine space
or whether the corresponding ring of invariant functions is a polynomial ring. Such quiver settings will be called coregular.

In this paper we present a method to determine if a random given
quiver setting $(Q, \alpha)$ is indeed coregular. Because the
quotient space $\iss_\alpha Q$ can be seen as the product of the
quotient spaces of the strongly connected components of
$(Q,\alpha)$ (see lemma \ref{connect}), we can restrict to
strongly connected quiver settings.

The method will consist of a number of allowed reduction steps.
Using these steps one attempts to simplify the quiver setting as
much as possible. When this is done one has to check whether the
reduced quiver setting is equal to one of $3$ basic quiver
settings that have a smooth quotient space. The main theorem we
will prove can be formulated as:

\begin{thm}\label{main}
Let $(Q,\alpha)$ be a genuine strongly connected quiver setting and $(Q',\alpha')$ is the quiver setting
obtained after all possible reductions of the form
\begin{enumerate}
\item[$\mathcal R_{I}$]
If $\sum_{j=1}^k i_j\le \alpha_v$ or $\sum_{j=1}^l u_j\le \alpha_v$ we delete the vertex $v$.
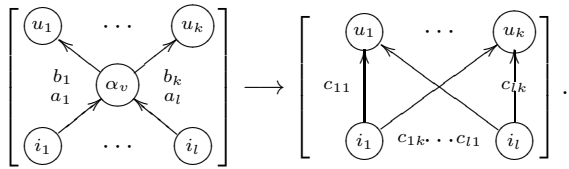
\[
\left[ \vcenter{
\xymatrix@=.3cm{
\vtx{u_1}&\cdots &\vtx{u_k}\\
&\vtx{\alpha_v}\ar[ul]^{b_1}\ar[ur]_{b_k}&\\
\vtx{i_1}\ar[ur]^{a_1}&\cdots &\vtx{i_l}\ar[ul]_{a_l}}}
\right]
\longrightarrow
\left[\vcenter{
\xymatrix@=.3cm{
\vtx{u_1}&\cdots &\vtx{u_k}\\
&&\\
\vtx{i_1}\ar[uu]^{c_{11}}\ar[uurr]_<<{c_{1k}}&\cdots &\vtx{i_l}\ar[uu]|{c_{lk}}\ar[uull]^<<{c_{l1}}}}
\right].
\]
\item[$\mathcal R_{II}$]
Remove the loops on a vertex with dimension $1$.
\[
\left[\vcenter{
\xymatrix@=.3cm{
&\vtx{1}\ar@{..}[r]\ar@{..}[l]\ar@(lu,ru)@{=>}^k&}}
\right]\longrightarrow
\left[\vcenter{
\xymatrix@=.3cm{
&\vtx{1}\ar@{..}[r]\ar@{..}[l]&}}
\right].\]
\item[$\mathcal R_{III}$]
Remove the only loop on a vertex with dimension $k>1$ which has a neighborhood like in one of the pictures below.
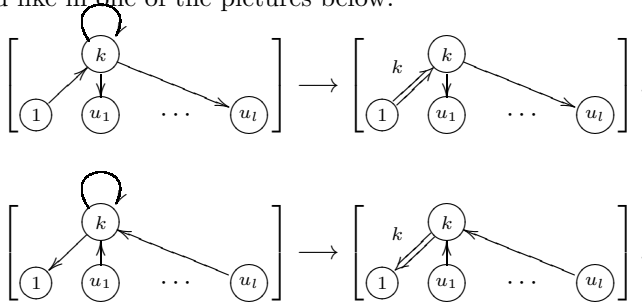
\[
\left[\vcenter{
\xymatrix@=.3cm{
&\vtx{k}\ar[d]\ar[drr]\ar@(lu,ru)&&\\
\vtx{1}\ar[ur]&\vtx{u_1}&\cdots &\vtx{u_l}}}
\right]\longrightarrow
\left[\vcenter{
\xymatrix@=.3cm{
&\vtx{k}\ar[d]\ar[drr]&&\\
\vtx{1}\ar@2[ur]^k&\vtx{u_1}&\cdots &\vtx{u_l}}}
\right],
\]
\vspace{.2cm}
\[
\left[\vcenter{
\xymatrix@=.3cm{
&\vtx{k}\ar[dl]\ar@(lu,ru)&&\\
\vtx{1}&\vtx{u_1}\ar[u]&\cdots &\vtx{u_l}\ar[ull]}}
\right]\longrightarrow
\left[\vcenter{
\xymatrix@=.3cm{
&\vtx{k}\ar@2[dl]_k&&\\
\vtx{1}&\vtx{u_1}\ar[u]&\cdots &\vtx{u_l}\ar[ull]}}
\right].
\]
\end{enumerate}
$(Q,\alpha)$ is coregular if and only if $(Q',\alpha')$ is one of the three settings below:
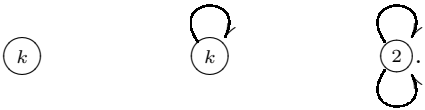
\[
\vspace{.5cm}
\xymatrix{\vtx{k}}\hspace{2cm} \xymatrix{\vtx{k}\ar@(lu,ru)}\hspace{2cm} \xymatrix{\vtx{2}\ar@(lu,ru)\ar@(ld,rd)}.
\]
\end{thm}

\section{Quiver representations}

In this section we recall some generalities about representations of quivers.
A \emph{quiver} $Q=(V,A,s,t)$ is a quadruple consisting of a set of vertices $V$, a set of arrows $A$  and 2
maps $s,t: A \to V$ which assign to  each arrow its starting and terminating vertex. We also denote this as
\begin{tiny}
\[
\xymatrix{\vtx{{t(a)}}&\vtx{{s(a)}}\ar[l]_a}.
\]
\end{tiny}
The \emph{Euler form} of $Q$ is the bilinear form
$\chi_Q : \Z^{\#V}\times \Z^{\#V} \to \Z$ defined by the matrix
\[
m_{ij} = \delta_{ij} -\#\{a|\xymatrix{\vtx{\scriptstyle {i}}&\vtx{\scriptstyle {j}}\ar[l]_a}
\},
\]
where $\delta$ is the Kronecker delta.
It is easy to see that that a quiver is uniquely defined by its Euler form.

A \emph{dimension vector} of a quiver is a map $\alpha: V \to \N$, the size of a dimension vector
is defined as $|\alpha| := \sum_{v\in V} \alpha_v$.
A couple $(Q,\alpha)$ consisting of a quiver and a dimension vector is called a \emph{quiver setting}
and for every vertex $v\in V$, $\alpha_v$ is refered to as the dimension of $v$.
If no vertex has dimension zero the setting is called \emph{genuine}.
If we draw pictures of quiver settings we will put the dimension of a vertex inside that vertex.

An \emph{$\alpha$-dimensional complex representation} $W$ of $Q$ assigns to each vertex
$v$ a linear space $\C^{\alpha_v}$ and to each arrow $a$ a matrix
\[
W_a \in \grp{Mat}_{\alpha_{t(a)} \times \alpha_{s(a)}}(\C).
\]
The space of all $\alpha$-dimensional representations is denoted by $\Repr_\alpha Q$.
\[
\Repr_\alpha Q := \bigoplus_{a \in A}\grp{Mat}_{\alpha_{t(a)} \times \alpha_{s(a)}}(\C).
\]
To the dimension vector $\alpha$ we can also assign a reductive group
\[
\GL_\alpha := \bigoplus_{v \in V}\grp{\GL_{\alpha_v}(\C)}.
\]
This group can be considered as the group of base changes in the vector spaces associated to the vertices.
Therefore every element of this group, $g$, has a natural action on $\Repr_\alpha Q$:
\[
W:=(W_a)_{a\in A},~W^g := (g_{t(a)}W_ag_{s(a)}^{-1})_{a\in A}
\]
Two representations in $\Rep_\alpha Q$ are called equivalent,
if they belong to the same orbit under the action of $\GL_\alpha$.

For every vertex we also define a special dimension vector
\[
\epsilon_v:V \to \N: w \mapsto \delta_{vw},
\]
and an $\epsilon_v$-dimensional representation $S_v$
assigning to every arrow the zero matrix.

A representation $W$ is called \emph{simple} if the only collections of subspaces
$(\mathsf V_v)_{v\in V},~\mathsf V_v \subseteq \C^{\alpha_v}$ having the property
\[\forall a\in A: W_a\mathsf V_{s(a)}\subset\mathsf V_{t(a)}\]
are the trivial ones (i.e. the collection of zero-dimensional subspaces
and ($\C^{\alpha_v})_{v \in V}$).

The direct sum $W\oplus W'$ of two representations $W, W'$ has as dimension vector the sum of the two dimension
vectors and as matrices $(W\oplus W')_a :=W_a\oplus W'_a$.
A representation equivalent to a direct sum of simple representations is called \emph{semisimple}.

From the algebraic point of view one can look at the ring of polynomial functions
over $\Rep_\alpha Q$ which is a polynomial ring denoted by $\C[\Repr_\alpha Q]$.
On this ring there is a corresponding action of $\GL_\alpha$ and one can look at the corresponding
subring of functions that are invariant under this action:
\[
\C[\Repr_\alpha Q]^{\GL_\alpha}:= \{f \in \C[\Repr_\alpha Q]|f^g=f\}.
\]
The variety corresponding to this subring is denoted by $\iss_\alpha Q$ and by \cite{Artin}
and \cite{Kraft:1984} this space classifies the equivalence classes
semisimple $\alpha$-dimensional representations of $Q$ which are in fact the closed
$\GL_\alpha$-orbits in $\Repr_\alpha Q$. The ring of invariants will also be
denoted by $\C[\iss_\alpha Q]$.

If $\iss_\alpha Q$ is a smooth variety then it is an affine space, this follows immediately from
(\cite{Kraft:1984} 4.3B lemma 1 p.139).
\begin{thm}\label{pol}
Suppose $V$ is a complex vector space with a linear action of a reductive group $\grp G$.
If the affine quotient $V/\grp G$ is smooth in the point corresponding to $0$ then
$V/G=\C^t$ for a $t \in \N$. The corresponding ring of invariants $\C[V]^\grp G$ is then a polynomial
ring.
\end{thm}

If we want to study the ring of invariants it is important to know by what functions it
is generated. The solution to this problem is given in the article by Le Bruyn and Procesi
about semisimple quiver representations \cite{LBP}.

A sequence of arrows $a_1\dots a_p$ in a quiver Q is called a \emph{path of length $p$}
if $s(a_i)=t(a_{i+1})$, this path is called a cycle if $s(a_p)=t(a_1)$.

To a cycle we can associate a polynomial function
\[
f_c:\Repr_{\alpha} Q \to \C : W \mapsto \grp{Tr}(W_{a_1}\cdots W_{a_p})
\]
which is definitely $\GL_{\alpha}$-invariant. Two cycles that are a cyclic permutation
of each other give the same polynomial invariant, because of the basic properties of the trace map.
Two such cycles are called equivalent.

A cycle $a_1\dots a_p$ is called \emph{primitive} if every arrow has a different starting vertex.
This means that the cycle runs through each vertex at most $1$ time. It is easy to see
that every cycle has a decomposition in primitive cycles. It is however not
true that the corresponding polynomial invariant decomposes to a product of the polynomial
functions of the primitive cycles.

We will call a cycle \emph{quasi-primitive} for a dimension vector $\alpha$ if the
vertices that are ran through more than once, have dimension bigger than $1$.
By cyclicly permuting a cycle and splitting the trace of a product of two $1\times 1$
matrices into a product of traces, we can always decompose an $f_c$ into a product of
traces of quasi-primitive cycles. We now have the following result

\begin{thm}[Le Bruyn-Procesi]\label{cycles}
$\C[\iss_\alpha Q]$ is generated by all $f_c$ where $c$ is a quasi-primitive cycle with length smaller than
$|\alpha|^2+1$.
We can turn $\C[\iss_\alpha Q]$ into a graded ring by giving $f_c$ the length
of its cycle as degree.
\end{thm}
This result can be used to prove and interesting lemma about the coregularity of
subquivers.
\begin{defn}
Define a partial ordering on the set of quivers in the following way.
A quiver $Q'=(V',A',s',t')$ is smaller than $Q=(V,A,s,t)$ if (up to isomorphism)
\[
V'\subseteq V,~A'\subseteq A,~s' =s|_{A'} \text{ and } t' =t|_{A'},
\]
$Q'$ is called a \emph{subquiver} of $Q$.
\end{defn}
\begin{lemma}\label{smoothsub}
If $\iss_\alpha Q$ is smooth and $Q'\le Q$ then
$\iss_{\alpha'} Q'$ is also smooth, where $\alpha':= \alpha|_{V'}$
\end{lemma}
\begin{prf}
We have an embedding
\[
\xymatrix{\Repr_{\alpha'} Q'\ar@{^{(}->}[r] & \Repr_{\alpha} Q}
\]
by assigning to the additional arrows in $Q$ zero matrices. So
\[
\xymatrix{\C[\Repr_{\alpha} Q]\ar@{->>}[r] & \C[\Repr_{\alpha'} Q']} \implies
\xymatrix{\C[\Repr_{\alpha} Q]^{\GL_\alpha}\ar@{->>}[r] & \C[\Repr_{\alpha'} Q']^{\GL_\alpha}}.
\]
Because the action of $\GL_\alpha$ on $\Repr_{\alpha'} Q'$ reduces to that of $\GL_{\alpha'}$,
$\C[\iss_{\alpha'} Q']$ is a quotient ring of $\C[\iss_\alpha Q]=\C[X_1,\dots,X_n]$.
The only relations that we have to divide out are the $X_i$ that correspond to
a cycle containing one of the additional arrows we put zero, so $\C[\iss_{\alpha'} Q']$
is just a polynomial ring with fewer variables.
\eop \end{prf}

Two vertices $v$ and $w$ are said to be \emph{strongly connected} if there is a path from
$v$ to $w$ and vice versa. It is easy to check that this relation is an equivalence
so we can divide the set of vertices into equivalence classes $V_i$. The subquiver $Q_i$
having $V_i$ as set of vertices, and as arrows all arrows between vertices of $V_i$
is called a \emph{strongly connected component} of $Q$.

\begin{lemma}\label{connect}
\begin{enumerate}
\item[]
\item
If $(Q,\alpha)$ is a quiver setting then
\[
\C[\iss_\alpha Q] := \bigotimes_{i} \C[\iss_{\alpha_i} Q_i]
\]
where $Q_i=(V_i,A_i,s_i,t_i)$ are the strongly connected components of $Q$ and $\alpha_i:=\alpha|_{V_i}$.
\item
$\iss_\alpha Q$ is smooth if and only if the $\iss_{\alpha} Q_i$ of all
its strongly connected components are smooth.
\end{enumerate}
\end{lemma}
\begin{prf}
\begin{enumerate}
\item[]
\item
By theorem \ref{cycles} $\C[\iss_\alpha Q]$ is generated by the traces of cycles. Every cycle
belongs to a certain connected component of $Q$. Between $f_c$'s coming from cycles of different components
there cannot be any relations, so we can consider the ring of invariants as a tensor-products of the rings of invariants
different strongly connected components.
\item
If all the strongly connected components are coregular the ring of invariants of the total quiver setting
will be the tensor product of polynomial rings and hence a polynomial ring. The inverse implication follows
directly from lemma \ref{smoothsub}.
\end{enumerate}
\eop \end{prf}

\section{Reduction Steps}

As we stated in the introduction we want to apply some kind of reduction
on quivers. By this we mean that if we start from a general quiver setting $(Q,\alpha)$,
we want to construct a new quiver setting with fewer vertices or arrows
but with the same or a closely related ring of invariants. In this section we will consider
3 different types of reductions.

First we have to recall a result from \cite{Kraft:1984}
\begin{thm}
Consider the vector space $\Mat_{k\times l}(\C) \oplus \Mat_{l\times m}(\C)$ together with an action
of $\GL_l(\C)$:
\[
(M_1,M_2)^g := (M_1g,g^{-1}M_2).
\]
The quotient space $\Mat_{k\times l}(\C) \oplus \Mat_{l\times m}(\C)/\GL_l(\C)$
is isomorphic to the space of all $k\times m$-matrices of rank smaller then $l$ (so if $l\ge k$ or $l \ge m$
there is no restriction on the matrices and the quotient space is $\Mat_{k\times m}(\C)$).
Identification happens via the $\GL_l(\C)$-invariant map
\[
\pi : (M_1,M_2) \mapsto M_1M_2.
\]
\end{thm}
This lemma can now be applied to quiver settings:

\begin{lemma}[Reduction $\mathcal R_I$: Removing Vertices]\label{reduction}
Suppose $(Q,\alpha)$ is a quiver setting and $v$ is a vertex without
loops such that
\[
\chi_Q(\alpha,\epsilon_v)\ge 0 \text{ or }\chi_Q(\epsilon_v,\alpha)\ge 0.
\]
Construct a new quiver setting $(Q',\alpha')$ by changing $Q$:
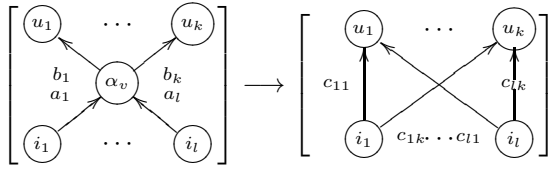
\begin{equation*}
\left[ \vcenter{
\xymatrix@=.3cm{
\vtx{u_1}&\cdots &\vtx{u_k}\\
&\vtx{\alpha_v}\ar[ul]^{b_1}\ar[ur]_{b_k}&\\
\vtx{i_1}\ar[ur]^{a_1}&\cdots &\vtx{i_l}\ar[ul]_{a_l}}}
\right]
\longrightarrow
\left[\vcenter{
\xymatrix@=.3cm{
\vtx{u_1}&\cdots &\vtx{u_k}\\
&&\\
\vtx{i_1}\ar[uu]^{c_{11}}\ar[uurr]_<<{c_{1k}}&\cdots
&\vtx{i_l}\ar[uu]|{c_{lk}}\ar[uull]^<<{c_{l1}}}} \right]
\end{equation*}
(some of the top and bottom vertices in the picture may be
the same).
These two quiver settings now have isomorphic rings of invariants.
\end{lemma}
\begin{prf}
We can split up the representation space into the following direct sum
\se{
\Repr_\alpha Q &=
\underbrace{\bigoplus_{a,~s(a)=v}\grp{Mat}_{\alpha_{t(a)} \times \alpha_{s(a)}}(\C)
}_{\text{arrows starting in }v}\oplus
\underbrace{\bigoplus_{a,~t(a)=v}\grp{Mat}_{\alpha_{t(a)} \times \alpha_{s(a)}}(\C)
}_{\text{arrows terminating in }v} \oplus
\text{ Rest}\\
&=
\grp{Mat}_{\sum_{s(a)=v}\alpha_{t(a)} \times \alpha_v}(\C)\oplus
\grp{Mat}_{\alpha_v \times \sum_{t(a)=v}\alpha_{s(a)} }(\C)\oplus
\text{ Rest}\\
&=\grp{Mat}_{\alpha_v-\chi(\alpha,\epsilon_v)\times \alpha_v}(\C)\oplus
\grp{Mat}_{\alpha_v \times \alpha_v-\chi(\epsilon_v,\alpha)}(\C)
\oplus
\text{ Rest}
}
The $\GL_{\alpha_v}(\C)$-part only acts on the first two terms and not on the
rest term. So if we take the quotient corresponding to $\GL_{\alpha_v}(\C)$
we only have to consider the first two terms.

By the previous lemma and keeping in mind that either $\chi_Q(\alpha,\epsilon_v)\ge 0$ or
$\chi_Q(\epsilon_v,\alpha)\ge 0$ the quotient space is equal to
\[
\grp{Mat}_{\alpha_v-\chi(\alpha,\epsilon_v)\times \alpha_v-\chi(\epsilon_v,\alpha)}(\C)\oplus
\text{ Rest}
\]
This space can be decomposed in the following way:
\[
\bigoplus_{\begin{matrix}a,&t(a)=v\\ b,&s(b)=v\end{matrix}}
\grp{Mat}_{\alpha_{t(b)}\times \alpha_{s(a)}}(\C) \oplus \text{ Rest}
\]
This direct sum is the same as the representation space
of the new quiver setting $(Q',\alpha')$.
\eop\end{prf}

\begin{lemma}[Reduction $\mathcal R_{II}$: Removing loops of dimension $1$]
Suppose that $(Q,\alpha)$ is a quiver setting and $v$  a vertex with $k$ loops and $\alpha_v=1$.
Take $Q'$ the corresponding quiver without loops, then the following identity hold
\[
\C[\iss_\alpha Q] \cong \C[\iss_\alpha Q'] \otimes \C[X_1,\cdots, X_k]
\]
\end{lemma}
\begin{prf}
This follows easily from $\ref{cycles}$ and the fact a cycle containing such a loop
can never be quasi-primitive unless it is the loop itself.
\eop\end{prf}

\begin{lemma}[Reduction $\mathcal R_{III}$: Removing a loop of higher dimension]
Suppose $(Q,\alpha)$ is a quiver setting and $v$ is a vertex of dimension $k \ge 2$ with
one loop such that
\[
\chi_Q(\alpha,\epsilon_v)= -1 \text{ or }\chi_Q(\epsilon_v,\alpha)=-1.
\]
Construct a new quiver setting $(Q',\alpha')$ by changing $(Q,\alpha)$:
\vspace{.5cm}
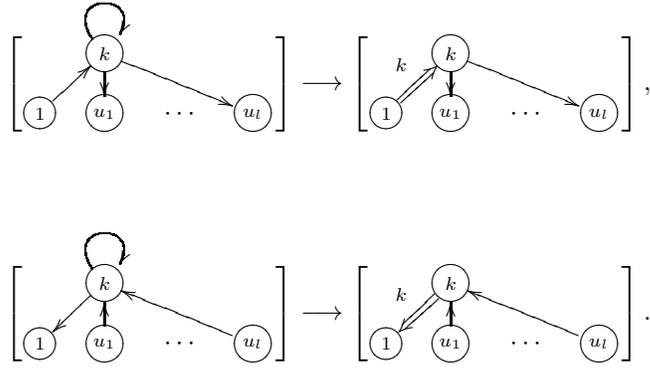
\[
\left[\vcenter{
\xymatrix@=.3cm{
&\vtx{k}\ar[d]\ar[drr]\ar@(lu,ru)&&\\
\vtx{1}\ar[ur]&\vtx{u_1}&\cdots &\vtx{u_l}}}
\right]\longrightarrow
\left[\vcenter{
\xymatrix@=.3cm{
&\vtx{k}\ar[d]\ar[drr]&&\\
\vtx{1}\ar@2[ur]^k&\vtx{u_1}&\cdots &\vtx{u_l}}}
\right],
\]
\vspace{.5cm}
\[
\left[\vcenter{
\xymatrix@=.3cm{
&\vtx{k}\ar[dl]\ar@(lu,ru)&&\\
\vtx{1}&\vtx{u_1}\ar[u]&\cdots &\vtx{u_l}\ar[ull]}}
\right]\longrightarrow
\left[\vcenter{
\xymatrix@=.3cm{
&\vtx{k}\ar@2[dl]_k&&\\
\vtx{1}&\vtx{u_1}\ar[u]&\cdots &\vtx{u_l}\ar[ull]}}
\right].
\]
We have the following identity:
\[
\C[\iss_\alpha Q]\cong \C[\iss_{\alpha'} Q']\otimes \C[X_1,\dots,X_k]
\]
\end{lemma}
\begin{prf}
We only prove this for the first case.
Call the loop in the first quiver $\ell$ and the incoming arrow $a$.  Call the incoming arrows in the second quiver $c_i, i=0,\dots,k-1$.

There is a map
\[
\pi: \Rep_\alpha Q \to \Rep_{\alpha'}Q'\times \C^k: V \mapsto
(V',\Tr V_\ell,\dots,\Tr V_{\ell}^k)\text{ with }V'_{c_i} :=
V_\ell^iV_a.
\]
Suppose $(V',x_1,\dots,x_k)\in \Rep_{\alpha'}Q'\times \C^k\in $
such that $(x_1,\dots,x_k)$ corresponds to the traces of powers of
an invertible diagonal matrix $D$ with $k$ different eigenvalues
($\lambda_i, i= 1,\dots,k$) and the matrix $A$ made of the columns
($V_{c_i}, i=0,\dots, k-1$) is invertible. The image of
representation
\[V \in \Rep_\alpha Q:
V_a = V'_{c_0},\vspace{.5cm}V_\ell =A
\left(
\begin{smallmatrix}
\lambda_1^{0}&\cdots&\lambda_1^{k-1}\\
\vdots&&\vdots\\
\lambda_k^{0}&\cdots&\lambda_k^{k-1}
\end{smallmatrix}
\right)^{-1}
D
\left(
\begin{smallmatrix}
\lambda_1^{0}&\cdots&\lambda_1^{k-1}\\
\vdots&&\vdots\\
\lambda_k^{0}&\cdots&\lambda_k^{k-1}
\end{smallmatrix}
\right)
A^{-1}
\]
under $\pi$ is $(V',x_1,\dots,x_k)$ because
\se{
V_\ell^iV_a &= A
\left(
\begin{smallmatrix}
\lambda_1^{0}&\cdots&\lambda_1^{k-1}\\
\vdots&&\vdots\\
\lambda_k^{0}&\cdots&\lambda_k^{k-1}
\end{smallmatrix}
\right)^{-1}
D^i
\left(
\begin{smallmatrix}
\lambda_1^{0}&\cdots&\lambda_1^{k-1}\\
\vdots&&\vdots\\
\lambda_k^{0}&\cdots&\lambda_k^{k-1}
\end{smallmatrix}
\right)
A^{-1}V'_{c_0}\\
&=A
\left(
\begin{smallmatrix}
\lambda_1^{0}&\cdots&\lambda_1^{k-1}\\
\vdots&&\vdots\\
\lambda_k^{0}&\cdots&\lambda_k^{k-1}
\end{smallmatrix}
\right)^{-1}
\left(
\begin{smallmatrix}
\lambda_1^{i}\\
\vdots\\
\lambda_k^{i}
\end{smallmatrix}
\right)\\
&=V_{c_i}
}
and the traces of $V_\ell$ are the same as the ones of $D$.
The conditions we imposed on $(V',x_1,\dots,x_k)$, imply that the image of $\pi$, $U$, is dense, and hence $\pi$ is a dominant map.

We have a bijection between the  generators of $\C[\iss_\alpha Q]$ and $\C[\iss_{\alpha'} Q']\otimes \C[X_1,\dots,X_k]$ by identifying
\[
f_{\ell^i} \mapsto X_i, i=1,\dots,k~ ,\vspace{.5cm} f_{\cdots a\ell^i\cdots}  \mapsto f_{\cdots c_i\cdots},i=0,\dots,k-1
\]
Notice that higher orders of $\ell$ don't occur because of
the Caley Hamilton identity on $V_\ell$.
So if $n$ is the number of generators of $\C[\iss_\alpha Q]$,
we have two maps
\se{
\phi &: \C[Y_1,\cdots Y_n] \to \C[\iss_\alpha Q] \subset \C[\Rep_\alpha Q],\\
\phi'&: \C[Y_1,\cdots Y_n] \to \C[\iss_{\alpha'} Q']\otimes \C[X_1,\dots,X_k]\subset \C[\Rep_{\alpha'} Q'\times \C^k].}

Notice that we have that
$\phi'(f) \circ \pi \equiv \phi(f)$ and
$\phi(f) \circ \pi^{-1}|_{U} \equiv \phi'(f)|_{U}$.
So if $\phi(f)=0$ then also $\phi'(f)|_{U}=0$. Because $U$ is zariski-open
and dense in $\Rep_{\alpha'} Q'\times \C^2$, $\phi'(f)\equiv 0$.
A similar argument holds for the inverse implication so $\Ker \phi = \Ker \phi'$.
\eop\end{prf}

We have seen three possible reductions of a quiver setting which
keep the ring of invariants intact or split of a tensor product
with a polynomial ring. We can also apply the inverse steps of the
reduction to add new vertices or loop while keeping the ring of
invariants the same or tensoring it up with a polynomial ring.
These inverse steps will be denoted as $\mathcal R_{\dots}^{-1}$.

The previous three lemma's can now be summarized as
\begin{thm}
Suppose that $(Q,\alpha)$  and $(Q',\alpha')$ are two quiver settings that can be transformed
into eachother using consecutive steps of the form $\mathcal R_{I}$, $\mathcal R_{I}^{-1}$,
$\mathcal R_{II}$, $\mathcal R_{II}^{-1}$, $\mathcal R_{III}$ or $\mathcal R_{III}^{-1}$.
Then $(Q,\alpha)$ is coregular if and only if $(Q',\alpha')$ is coregular.
\end{thm}

\begin{defn}
A quiver setting $(Q,\alpha)$ such that there cannot be applied any reduction steps
$\mathcal R_{I}$, $\mathcal R_{II}$ or $\mathcal R_{III}$ will be called \emph{reduced}.
\end{defn}

It remains now to search for the reduced coregular quiver settings.
As we will see there are only a very limited number of them.
But before we do that we must introduce some techniques that allow us to
rule out non coregular quiver settings.

\section{Local Quiver settings}
The technique of local quiver settings is very useful to rule out quiver settings that
are not coregular.
If we want to prove that a certain  $(Q,\alpha)$ is coregular,
we have to check that $\iss_\alpha Q$ is smooth in every point.
Take a point $p \in \iss_\alpha Q$, this point will correspond to
the isomorphism class of a semisimple representation $V \in \Repr_\alpha Q$ which can be decomposed as a direct sum of
simple representations.
\[
V= S_1^{\oplus a_1}\oplus \dots \oplus S_k^{\oplus a_k},
\]
A theorem by Le Bruyn and Procesi \cite[Theorem 5]{LBP} states that we can build a new quiver setting with a similar quotient space, but having a simpler
structure.

\begin{thm}[Le Bruyn-Procesi]\label{local}
For a point $p \in \iss_\alpha Q$ corresponding to a semisimple representation $V=S_{1}^{\oplus a_1}\oplus \dots \oplus S_{k}^{\oplus a_k}$,
there is a quiver setting $(Q_p,\alpha_p)$ called the \emph{local quiver setting} such that
we have an \'etale isomophism between an open
neighborhood of the zero representation in $\iss_{\alpha_p} Q_p$ and an open neighborhood of $p$.

$Q_p$ has $k$ vertices corresponding to the set $\{S_i\}$ of simple factors of $V$ and between
$S_i$ and $S_j$ the number of arrows equals
\[
\delta_{ij}-\chi_Q(\alpha_i,\alpha_j)
\]
where $\alpha_i$ is the dimension vector of the simple component $S_i$ and $\chi_Q$ is the Euler form
of the quiver $Q$.
The dimension vector $\alpha_p$ is defined to be $(a_1,\dots,a_k)$, where the $a_i$ are the multiplicities of the simple components in $V$.
\end{thm}

Suppose now that we want to find out whether a certain
space $\iss_\alpha Q$ is smooth. If this were the case we can choose a certain point $p$
and look at it locally. Because of the \'etale isomorphism,
the corresponding local quiver $Q_p$ must
have a quotient space $\iss_{\alpha_p}Q_p$ that is smooth in the zero representation.
Therefore by \ref{pol}, $\C[\iss_{\alpha_p}Q_p]$ must be a polynomial ring
and hence $(Q_p,\alpha_p)$ is coregular. This must hold for every point so we have to
check all possible points $p$.

\begin{thm}\label{smoothlocal}
$(Q,\alpha)$ is coregular if and only if for every possible semisimple $\alpha$-dimensional
representation $V$, the corresponding
local quiver setting is coregular.
\end{thm}

One of the local quivers is equal to the original quiver, namely the one corresponding
to the $\alpha$-dimensional zero-representation
\[
\bigoplus_{v \in V} S_v^{\oplus\alpha_v},
\]
This implies that we can only use this result to rule out quiver settings that are not coregular.

The structure of the local quiver setting only depends on the dimension vectors
of the simple components. Therefore one can restrict to looking at decompositions
of $\alpha$ into dimension vectors $\beta_i$ f.i.
\[
\alpha = a_1\beta_1+\cdots + a_k\beta_k\text{ (the $\beta_i$ need not to be different).}
\]
One can now ask whether there is a semisimple representation
corresponding to such a decomposition.
The answer to this question will be positive whenever for all
the $\beta_i$ there exist simple representations of that dimension vector and
if there are two or more $\beta_i$ equal, there are at least as many different simple
representation classes with dimension vector $\beta_i$ (otherwise you cannot make a
direct sum with different simple representations having the same dimension vector).

To check the above conditions we must also have a characterization of the dimension vectors for
which a quiver has simple representations. We recall a result from Le Bruyn
and Procesi \cite[Theorem 4] {LBP}.
\begin{thm}\label{simple}
Let $(Q, \alpha)$ be a genuine quiver setting.
There exist simple representations of dimension vector $\alpha$ if and only if
\begin{itemize}
\item
If $Q$ is of the form
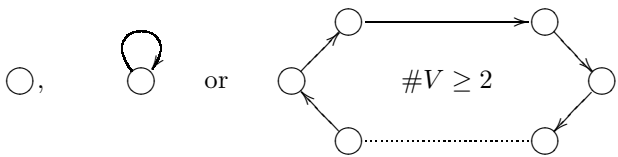
\[
\vcenter{\xymatrix@=.5cm{\vtx{~}}},\hspace{1cm}
\vcenter{\xymatrix@=.5cm{\vtx{~}\ar@(lu,ru)}}\hspace{.5cm}\text{ or }\hspace{.5cm}
\vcenter{\xymatrix@=.3cm{
&\vtx{~}\ar@{->}[rr] &&\vtx{~}\ar@{->}[rd] &\\
\vtx{~}\ar@{->}[ru] &&\#V\ge 2&&\vtx{~}\ar@{->}[ld] \\
&\vtx{~}\ar@{->}[lu] &&\vtx{~}\ar@{..}[ll]&\\
}}
\]
and $\alpha = \unit$ (this is the constant map from the vertices to $1$).
\item
$Q$ is not of the form above, but strongly connected and
\[
\forall v \in V:\chi_Q(\alpha,\epsilon_v)\le 0 \text{ and }\chi_Q(\epsilon_v,\alpha)\le 0
\]
(we recall that a quiver is straongly connected if and only if between every two vertices there
are paths connection them in both directions).
\end{itemize}
In both cases the dimension of $\iss_\alpha Q$ is given by $1-\chi_Q(\alpha,\alpha)$.
In all cases except for the one vertex without loops this dimension is bigger then $0$, so
then there are infinite classes of simples with that dimension vector.
In the case of the one vertex $v$ without loops, there is one unique simple representation
$S_v$.

If $(Q,\alpha)$ is not genuine, the simple representations classes
are in bijective correspondence to the simple representations classes
of the genuine quiver setting obtained by deleting all vertices with dimension zero.
\end{thm}

To rule out quiver settings that are not coregular we must find a local quiver
setting that is not coregular or contains a non-coregular subquiver setting
by lemma \ref{smoothsub}.

For symmetric quiver settings, these are quiver settings
with a symmetric Euler form, \cite{raf} gives us a complete classification
of all possible quiver settings that are coregular.

\begin{defn}
A quiver $Q=(V,A,s,t)$ is said to be the \emph{connected sum} of $2$ subquivers $Q_1=(V_1,A_1,s_1,t_1)$ and $Q_1=(V_2,A_2,s_2,t_2)$ at the vertex $v$,
if the two subquivers make up the whole quiver and only intersect in the vertex $v$.
So in symbols $V= V_1 \cup V_2$, $A= A_1 \cup A_2$, $V_1 \cap V_2=\{v\}$ and $A_1 \cap A_2=\emptyset$.
\[
Q_1 \substack{{\#}\\\scriptstyle{v}} Q_2:=
\vcenter{\xymatrix {
\dots\ar[rd]\ar@{..}@/^17mm/[dd]&&\dots\ar[ld]\ar@{..}@/_17mm/[dd]\\
Q_1&\vtx {v}\ar[rd]\ar[ld]&Q_2\\
\dots&&\dots}}
\]
If we connect three or more components we write $Q_1 \substack{{\#}\\\scriptstyle{v}} Q_2 \substack{{\#}\\\scriptstyle{w}} Q_3$ instead of
$(Q_1 \substack{{\#}\\\scriptstyle{v}} Q_2) \substack{{\#}\\\scriptstyle{w}} Q_3$ for sake of simplicity.
\end{defn}

\begin{thm}\label{symm}
Let $(Q,\alpha)$ be a symmetric strongly connected quiver setting without.
Then $(Q,\alpha)$ is coregular if and only if
$Q$ is a connected sum
\[
Q := Q_1 \substack{{\#}\\\scriptstyle{v_1}} Q_2 \substack{{\#}\\\scriptstyle{v_2}} \cdots \substack{{\#}\\\scriptstyle{v_{l-1}}} Q_l,
\]
where the $(Q_i,\alpha_i)$ are of the form
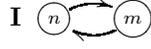
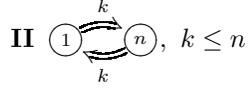
\begin{itemize}
\item[{\bf I}]
$\xymatrix@=.5cm{
\vtx{n}\ar@/^/[r]&\vtx{m}\ar@/^/[l]
}$\\
\item [{\bf II}]
$\xymatrix@=.5cm{
\vtx{1}\ar@2@/^/[r]^k&\vtx{n}\ar@2@/^/[l]^k
}, ~k\le n$\\
\item[{\bf III}]
$\xymatrix@=.5cm{
\vtx{1}\ar@/^/[r]&\vtx{n}\ar@/^/[l]\ar@/^/[r]&\vtx{m}\ar@/^/[l]
}$\\
\item[{\bf IV}]
$\xymatrix@=.5cm{
\vtx{n}\ar@/^/[r]&\vtx{2}\ar@/^/[l]\ar@/^/[r]&\vtx{m}\ar@/^/[l]
},
$\\
\end{itemize}
and $\alpha_{v_j}=1,~j=1,\dots,l-1$
\end{thm}

\section{Reduced coregular quiver settings}

First we look a the case of loops
\begin{lemma}\label{loops}
Suppose $(Q,\alpha)$ is a coregular strongly connected quiver setting such that
\[
\forall w \in V:\chi_Q(\alpha,\epsilon_w)< 0\text{ and }\chi_Q(\epsilon_w,\alpha)< 0.
\]
If $v$ is a vertex with loops
then $\alpha_v= 1$ or the neighborhood of $v$ has the following form
\vspace{.5cm}
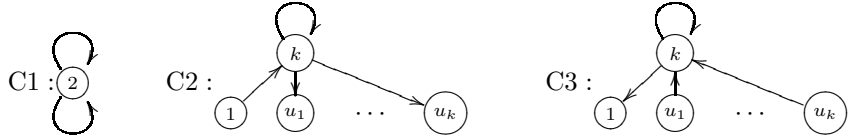
\[
{\rm C1:}
\vcenter{
\xymatrix@=.3cm{
\vtx{2}\ar@(lu,ru)\ar@(ld,rd) }}
\text{\hspace{1cm}}
{\rm C2:}
\vcenter{
\xymatrix@=.3cm{
&\vtx{k}\ar[d]\ar[drr]\ar@(lu,ru)&&\\
\vtx{1}\ar[ur]&\vtx{u_1}&\cdots &\vtx{u_k}}}
\text{\hspace{1cm}}
{\rm C3:}
\vcenter{
\xymatrix@=.3cm{
&\vtx{k}\ar[dl]\ar@(lu,ru)&&\\
\vtx{1}&\vtx{u_1}\ar[u]&\cdots &\vtx{u_k}\ar[ull]}}
\]
\vspace{.5cm}
\end{lemma}
\begin{prf}
\fbox{1. if $\alpha_v\ge 3$ there is only one loop in $v$}\\
Suppose that $\alpha_v\ge 3$ there are at least two loops in $v$.
In this case we have a subquiver as shown below. This subquiver
can be transformed into a symmetric quiver without loops using
lemma \ref{reduction} (in both ways). By \ref{symm} this symmetric
setting is not coregular, if $\alpha_v>2$. \vspace{.5cm}
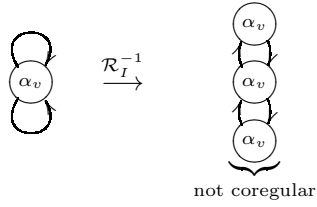
\[
\vcenter{
\xymatrix@=.3cm{
\vtx{\alpha_v}\ar@(lu,ru)\ar@(ld,rd)}}
\text{\hspace{.5cm}}\stackrel{{\cR_{I}^{-1}}}{\longrightarrow}
\text{\hspace{.5cm}}
\underbrace{\vcenter{
\xymatrix@=.3cm{
\vtx{\alpha_v}\ar@/^/[d]\\
\vtx{\alpha_v}\ar@/^/[u]\ar@/^/[d]\\
\vtx{\alpha_v}\ar@/^/[u]}}}_{\text{not coregular}}
\]

\fbox{2. If $\alpha_v= 2$ we are in C1 or there is only 1 loop in $v$}\\

If $\alpha_v= 2$ and we or not in $C1$, $C2$ or $C3$, $Q$ has either at least $3$ loops or either two loops and a cyclic path through $v$ (this cyclic path
can be constructed because $Q$ is strongly connected and contains at least $2$ vertices, otherwise $(Q,\alpha)=C1$).

In both cases we can take again the corresponding subquivers and change them to a symmetric quiver without loops
which is not coregular according to \ref{symm}.
\vspace{.5cm}
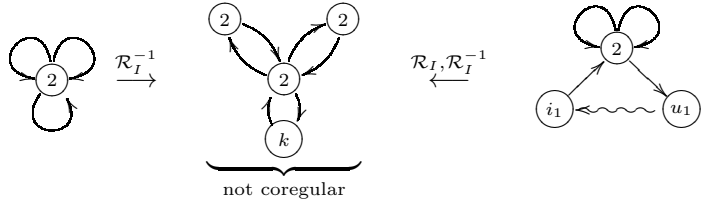
\[
\vcenter{
\xymatrix@=.3cm{
\vtx{2}\ar@(u,r)\ar@(u,l)\ar@(dl,dr) }}
\text{\hspace{.5cm}}\stackrel{{\cR_{I}^{-1}}}{\longrightarrow}
\text{\hspace{.5cm}}
\underbrace{\vcenter{
\xymatrix@=.3cm{
\vtx{2}\ar@/^/[dr]&&\vtx{2}\ar@/^/[dl]\\
&\vtx{2}\ar@/^/[ur]\ar@/^/[ul]\ar@/^/[d]&\\
&\vtx{k}\ar@/^/[u]&}}}_{\text{not coregular}}
\text{\hspace{.5cm}}\stackrel{{\cR_{I},\cR_{I}^{-1}}}{\longleftarrow}
\text{\hspace{.5cm}}
\vcenter{
\xymatrix@=.3cm{
&\vtx{2}\ar[dr]\ar@(u,r)\ar@(u,l)&&\\
\vtx{i_1}\ar[ur]& &\vtx{u_1}\ar@{~>}[ll]}}
\]
So the only possibility with more than one loop is C1.

\fbox{3. If $\alpha_v\ge 2$ and there is only 1 loop in $v$ then we are in C2 or C3}\\

Suppose that the dimension in $v$ is bigger than $1$ and that there is only $1$ loop.
Consider the representation
\[
W \oplus L \oplus \left(\bigoplus_{w\in V}S_w^{\oplus \alpha_w-1-\delta_{vw}}\right)
\]
where $W$ is a simple representation with dimension vector $\unit$
which is the constant map assigning $1$ to every vertex. Such a
representation exists by \ref{simple} because $Q$ is strongly
connected and $\chi_Q(1,\epsilon_w)\le 0$. $S_w$ is the
representation with dimension vector $\epsilon_w$ which assigns to
every arrow a zero matrix, while $L$ is a representation with
dimension vector $\epsilon_v$ which assigns to the loop in $v$ a
non-zero matrix.

For every vertex $w\ne v$ with dimension bigger than $1$ the local quiver contains exactly one vertex corresponding to the simple
representation $S_w$. For $v$ there is at least one vertex in the local quiver coming from $L$, which has dimension $1$.
If $\alpha_v>2$ there is an extra vertex from the $S_v$ but we won't consider it because it doesn't change the proof.

The subquiver containing the vertices from $L$ en $S_w,w\ne v$ is the same as in the original quiver
because
\[
\chi_Q(\epsilon_u,\epsilon_w) = \delta_{uw} - \#\{a|\xymatrix{\vtx{\scriptstyle {u}}&\vtx{\scriptstyle {w}}\ar[l]_a}\}
\]

In the local quiver we will draw the additional vertex coming from $W$ as a square.
The number of arrows from another vertex coming from $S_w$ to
the vertex coming from $W$ is equal to $-\chi_Q(\unit,\epsilon_w)$
and hence one less than the number of arrows
leaving $w$ in the original quiver.
The same holds for the number of arrows in the opposite direction and for the arrows between $L$ and $W$.

We will now look closely at the neighborhood of $v$.
\begin{itemize}
\item
\fbox{$\chi_Q(\epsilon_v,\unit)\le -2$ and $\chi_Q(\unit,\epsilon_v)\le -2$ is impossible}\\

The local quiver has a subquiver containing $\xymatrix{\vtx{1}\ar@2@/^/[r]&\vierkant{1}\ar@2@/^/[l]}$,
and $(Q,\alpha)$ is not coregular. For $(Q,\alpha)$ to be a coregular quiver setting,
one can suppose that either $\chi_Q(\epsilon_v,\unit)=-1$ or $\chi_Q(\unit, \epsilon_v)=-1$.
\vspace{.25cm}

\item
\fbox {$\chi_Q(\epsilon_v,\unit)=-1$ and $\chi_Q(\unit,\epsilon_v)\le -2$ implies C2.}\\

We claim that if $w_1$ is the unique vertex in $Q$ such that $\chi_Q(\epsilon_v,\epsilon_{w_1})=-1$ then
$\alpha_{w_1}=1$.

If this was not the case there is a vertex correponding to $S_{w_1}$ in the local quiver.
If $\chi_Q(\unit ,\epsilon_{w_1})=0$ then the dimension of the unique vertex $w_2$ with an arrow to $w_1$
has strictly bigger dimension than $w_1$, otherwise $\chi_Q(\alpha,\epsilon_{w_1})\ge 0$. The vertex $w_2$ corresponds again to
a vertex in the local quiver. If $\chi_Q(\unit ,\epsilon_{w_2})=0$, the unique vertex $w_3$ with an arrow to $w_2$
has strictly bigger dimension than $w_2$.
Proceeding this way one can find a sequence of vertices with increasing dimension, which attains a maximum
in vertex $w_k$. Therefore $\chi_Q(\unit ,\epsilon_{w_k})\le -1$.
This last vertex is in the local quiver connected with $W$,
so one has a path from $\unit$ to $\epsilon_v$.
\vspace{.5cm}
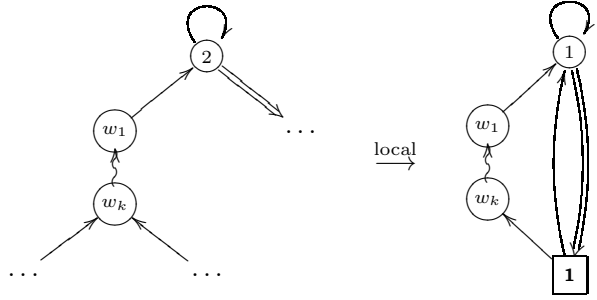
\[
\vcenter{\xymatrix@=.5cm{
&&\vtx{2}\ar@(lu,ru)\ar@{=>}[rd]&\\
&\vtx{w_1}\ar[ur]&&\dots\\
&\vtx{w_k}\ar@{~>}[u]&&\\
\dots\ar[ur]&&\dots\ar[ul]&}}
\text{\hspace{.5cm}}\stackrel{\text{local}}{\longrightarrow}
\text{\hspace{.5cm}}
\vcenter{\xymatrix@=.5cm{
&\vtx{1}\ar@(lu,ru)\ar@/^/@{=>}[ddd]&\\
\vtx{w_1}\ar[ur]&&\\
\vtx{w_k}\ar@{~>}[u]&&\\
&\vierkant{\unit}\ar[ul]\ar@/^/[uuu]&}}
\]
The local subquiver consisting of the vertices corresponding to $W$, $S_v$ and the $S_{w_i}$ is reducible via $\mathcal R_I$ to
$\xymatrix{\vtx{1}\ar@2@/^/[r]&\vierkant{1}\ar@2@/^/[l]}$. So if $\alpha_{w_1}>1$, $(Q,\alpha)$ is not coregular.
\vspace{.25cm}

\item
\fbox{$\chi_Q(\epsilon_v,\unit)\le -2$ and $\chi_Q(\unit,\epsilon_v)= -1$ implies C3.}\\

This follows by symmetry.\vspace{.25cm}

\item
\fbox{$\chi_Q(\epsilon_v,\unit)=-1$ and $\chi_Q(\epsilon_v,\unit)= -1$ implies C2 or C3.}\\
Suppose $w_1$ is the unique vertex in $Q$ such that $\chi_Q(\epsilon_v,\epsilon_{w_1})=-1$
and $w_k$ is the unique vertex in $Q$ such that $\chi_Q(\epsilon_{w_k},\epsilon_v)=-1$,
then either $\alpha_{w_1}=1$ or $\alpha_{w_k}=1$.

If this was not the case, consider the path connecting $w_k$ and
$w_1$ and call the intermediate vertices $w_i,~1<i<k$. Starting
from $w_1$ we go back along the path until $\alpha_{w_i}$ reaches
a maximum. At that point we know that $\chi_Q(\unit
,\epsilon_{w_k})\le -1$, otherwise $\chi_Q(\alpha
,\epsilon_{w_k})\ge 0$. In the local quiver there is a path from
the vertex corresponding to $W$ over the ones from $S_{w_i}$ to
$S_v$. Doing the same thing starting from $w_k$ we also have a
path from the vertex from $S_v$ over the ones of $S_{w_j}$ to $W$.

\vspace{.5cm}
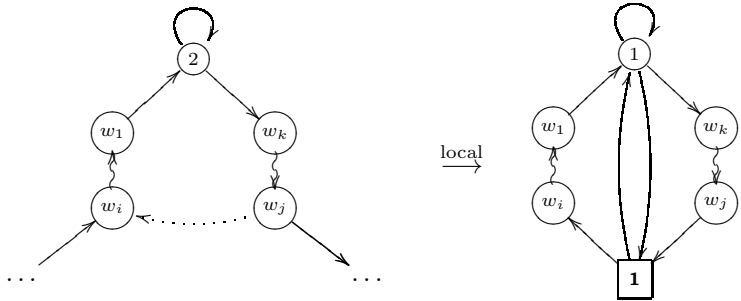
\[
\vcenter{\xymatrix@=.5cm{
&&\vtx{2}\ar@(lu,ru)\ar[rd]&&\\
&\vtx{w_1}\ar[ur]&&\vtx{w_k}\ar@{~>}[d]&\\
&\vtx{w_i}\ar@{~>}[u]&&\vtx{w_j}\ar[dr]\ar@/^/@{..>}[ll]\ar[dr]&\\
\dots\ar[ur]&&&&\dots}}
\text{\hspace{.5cm}}\stackrel{\text{local}}{\longrightarrow}
\text{\hspace{.5cm}}
\vcenter{\xymatrix@=.5cm{
&\vtx{1}\ar@(lu,ru)\ar@/^/[ddd]\ar[dr]&\\
\vtx{w_1}\ar[ur]&&\vtx{w_k}\ar@{~>}[d]\\
\vtx{w_i}\ar@{~>}[u]&&\vtx{w_j}\ar[dl]\\
&\vierkant{\unit}\ar[ul]\ar@/^/[uuu]&}}
\]
The subquiver consisting of $\unit$, $\epsilon_v$ and the two
paths through the $\epsilon_{w_i}$ is reducible to
$\xymatrix{\vtx{1}\ar@2@/^/[r]&\vierkant{1}\ar@2@/^/[l]}$. So if
both $\alpha_{w_1}>1$ and $\alpha_{w_k}>1$, $(Q,\alpha)$ is not
coregular.
\end{itemize}
\eop\end{prf}

We will now look at the reduced quiver settings without loops.

\begin{lemma}
A quiver setting with dimension vector $\unit$ is coregular if and only if the number of
primitive cycles equals the dimension of $\C[\iss_\unit Q]$.
\end{lemma}
\begin{prf}
The condition is obviously sufficient. It is also necessary
because if the number of cycles is bigger than the dimension then
there will be a relation between the cycles. If $\C[\iss_\unit Q]$
is a polynomial ring, these relations must be of the form
$Y=X_1\dots X_k$ but this is impossible because $Y$ is a primitive
cycle. \eop
\end{prf}

\begin{lemma}\label{noncoreg}
A strongly connected reduced quiver setting without loops is never coregular.
\end{lemma}
\begin{prf}
If $\alpha \ne \unit$,
consider the vertex $v$ with the highest
dimension. Then there exists indeed simple representations with dimension vector $\alpha-\epsilon_v$
because a reduced setting is never of the form
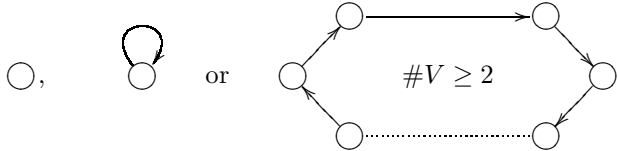
\[
\vcenter{\xymatrix@=.5cm{\vtx{~}}},\hspace{1cm}
\vcenter{\xymatrix@=.5cm{\vtx{~}\ar@(lu,ru)}}\hspace{.5cm}\text{ or }\hspace{.5cm}
\vcenter{\xymatrix@=.3cm{
&\vtx{~}\ar@{->}[rr] &&\vtx{~}\ar@{->}[rd] &\\
\vtx{~}\ar@{->}[ru] &&\#V\ge 2&&\vtx{~}\ar@{->}[ld] \\
&\vtx{~}\ar@{->}[lu] &&\vtx{~}\ar@{..}[ll]&\\
}}
\]
and $\alpha-\epsilon_v$ satisfies the second condition of theorem \ref{simple}:
\begin{itemize}
\item
If there is no arrow from $w$ to $v$,
$\chi_Q(\alpha-\epsilon_v,\epsilon_w)=\chi_Q(\alpha,\epsilon_w)\le
-1$.
\item
If there are $k$ arrows from $w$ to $v$ then $\chi_Q(\alpha,\epsilon_w)\le \alpha_w-k\alpha_v\le (1-k)\alpha_v $ so
\[
\chi_Q(\alpha-\epsilon_v,\epsilon_w)\le (1-k)\alpha_v +\chi_Q(\epsilon_v,\epsilon_w)= (1-k)\alpha_v -k \le -1.
\]
\item
Finally for $v=w$
\[
\chi_Q(\alpha-\epsilon_v,\epsilon_v)=\chi(\alpha,\epsilon_v)-1<-1\text{
and }\chi_Q(\epsilon_v,\alpha-\epsilon_v)\le -1.
\]
\end{itemize}
For reasons of symmetry $\chi_Q(\epsilon,\alpha-\epsilon_v)$ will
also be smaller than $0$ for every $w\in V$.

Due to the inequality $\chi_Q(\epsilon_v,\alpha-\epsilon_v)\le
-1$, the local quiver of a decomposition of the form
\[
(Q,\alpha-\epsilon_v) \oplus (Q,\epsilon_v)
\]
will not be coregular.

Suppose thus $\alpha=\unit$.
Because $(Q,\alpha)$ is reduced, there are at least $2$ arrows arriving and leaving every vertex.
For a connected quiver without loops $\Dim \C[\iss_\unit Q]=\#A-\#V +1$ so we have to prove
that for such quivers the number of primitive cycles is bigger than $\#A-\#V +1$ or that $Q$
constains a subquiver that is not coregular. We will do this by induction on the vertices.

\begin{itemize}
\item
For $\#V=2$ the statement is true because
\[Q := \xymatrix{\vtx{1}\ar@2@/^/[r]^k&\vtx{1}\ar@2@/^/[l]^l},~k,l\ge 2 ~\implies~ kl>k+l-1 .\]

\item
Suppose $\#V>2$ and that we have a subquiver of the form
\[\xymatrix{\vtx{1}\ar@2@/^/[r]^k&\vtx{1}\ar@2@/^/[l]^l}~(*)\]
If $k,l>1$ we know that this subquiver is not coregular and hence
neither is $Q$.

If both $k$ and $l$ are $1$ then replace this subquiver by $1$ vertex.
\[
\left[
\vcenter{\xymatrix@=.3cm{&&&\\
\vdots&\vtx{1}\ar@/^/[r]\ar[ul]&\vtx{1}\ar@/^/[l]\ar[ur]&\vdots\\
{}\ar[ur]&&&\ar[ul]
}}
\right]
\longrightarrow
\left[
\vcenter{\xymatrix@=.3cm{&&\\
\vdots&\vtx{1}\ar[ul]\ar[ur]&\vdots\\
{}\ar[ur]&&\ar[ul]
}}
\right]
\]
The new quiver $Q'$ is again reduced without loops because there
are at least $4$ arrows arriving in one of the vertices of the
subquiver and we only deleted $2$, the same holds for the arrows
leaving the subquiver. $Q'$ has one primitive cycle less than the
original. By induction we have that \se{
\Dim \C[\iss_\unit Q]&=\Dim \C[\iss_\unit Q']+1\\
&>(\#A' -\#V' +1) +1\\
&= \#A -\#V +1.
}
If for instance $k>1$ then one can look at the subquiver of $Q$ obtained by deleting
the $k-1$ edges, if this quiver is reduced then we are in the previous situation.
If this is not the case $Q$ contains a subquiver of the form
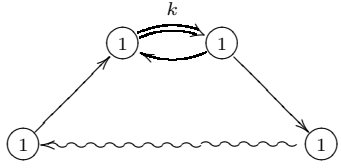
\[
\xymatrix{&\vtx{1}\ar@2@/^/[r]^k&\vtx{1}\ar@/^/[l]\ar[dr]&\\
\vtx{1}\ar[ur]&&&\vtx{1}\ar@{~>}[lll]
},
\]
which is not coregular because it is reducible to $(*)$.

\item
If $\#V>2$ and there are no subquivers of the form $(*)$, we can consider an arbitrary vertex $v$.
Construct a new quiver $Q'$ by performing the following substitution for $v$
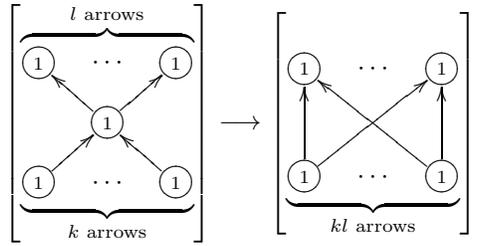
\[\left[
\overbrace{
\underbrace{ \vcenter{
\xymatrix@=.3cm{
\vtx{1}&\cdots &\vtx{1}\\
&\vtx{1}\ar[ul]\ar[ur]&\\
\vtx{1}\ar[ur]&\cdots &\vtx{1}\ar[ul]}}
}_{k \text{ arrows}}}^{l \text{ arrows}}
\right]
\longrightarrow
\left[\underbrace{
\vcenter{
\xymatrix@=.3cm{
\vtx{1}&\cdots &\vtx{1}\\
&&\\
\vtx{1}\ar[uu]\ar[uurr]&\cdots &\vtx{1}\ar[uu]\ar[uull]}}}
_{kl \text{ arrows}}\right].
\]
$Q'$ is again reduced without loops and
has the same number of primitive cycles, so by induction
\se{
\Dim \C[\iss_\unit Q]&=\Dim \C[\iss_\unit Q']\\
&>\#A' -\#V' +1 \\
&= \#A+(kl-k-l) -\#V+1 +1\\
&> \#A -\#V+1.
}
\end{itemize}
\eop\end{prf}

All this leads to the proof of our main theorem.
\begin{prf}
Statement \ref{main} follows immediately from lemmas \ref{loops}
and \ref{noncoreg} and the fact that as proven in \cite{Procesi}
the quiver settings that are listed in the theorem are coregular
\eop \end{prf} \emph{Acknowledgment.} I thank N. Verlinden for his
generous help in the graph theoretical part of proof
\ref{noncoreg}. \nocite{ALB} \nocite{King94}

\end{document}